\documentclass[12pt]{amsart}
\usepackage{amsfonts}
\usepackage{amssymb}

\usepackage{amsmath}
\usepackage{epsf}
\usepackage{epsfig}
\usepackage{colordvi}
\usepackage{fancybox}

\newtheorem{thm}{Theorem}[section]
\newtheorem{cor}[thm]{Corollary}

\theoremstyle{definition}

\theoremstyle{remark}

\newcommand{\ds}{\displaystyle}

\begin{document}

\title[Congruence of minimal surfaces  ]
{Congruence of minimal surfaces}%

\thanks
{This research is supported by the National Science Fund 
of the Bulgarian Ministry of Education and Science 
under Grant~DFNI-T01/0001, 2012.}
\author{Ognian Kassabov}%
\address{University of Transport, Sofia, Bulgaria}
\email{okassabov@abv.bg}


\begin{abstract}
An interesting problem in classical differential geometry
is to find methods to prove that two surfaces defined by different 
charts actually coincide up to position in  space. In a previous 
paper we proposed a method in this direction for minimal surfaces. 
Here we explain  not only how this method works but also how we can 
find the correspondence between the minimal surfaces, if they are 
congruent.  We show that two  families of minimal surfaces
which are proved to be conjugate actually coincide and  coincide
with their associated surfaces. We also consider  another family
of minimal polynomial surfaces of degree 6 and we apply the method
to show that some of them are congruent.
\end{abstract}
\maketitle

\section{Introduction}

Minimal surfaces have a long and interesting history. The name {\it minimal}
corresponds to the property that each point of such a surface 
has a neighborhood that minimizes the area among all surfaces with the same
boundary.  It turns out that this is equivalent to the vanishing of 
the mean curvature of the surface. The first nontrivial minimal surfaces discovered 
were the catenoid and the helicoid. After that, for many years no new minimal surfaces
were found. From  the investigations of Weierstrass and Enneper in the 19th century 
we know that any minimal surface can be found at least locally from a 
complex minimal curve, see e.g. \cite{Eis}. 

Nowadays the minimal surfaces are important in many areas such as 
architecture, material science, aviation,  biology. Because of their 
minimizing property, these surfaces can be also interesting  in CAD
research. Note that, for the moment, in this area it is more important to 
consider polynomial minimal surfaces. The problem is that such surfaces
are probably not very numerous in small degrees. In \cite{C-M} it is shown that
up to homothety and position in the space,  the classical Enneper surface is
the unique polynomial minimal surface of degree 3. 

In  \cite{X-W-1} and \cite{X-W-2} some families of polynomial minimal surfaces 
of degrees 6 and 5 are introduced. In \cite{OK} we showed a relation among these 
families. To do this, we cannot use the classical method used in \cite{C-M}
for degree 3. We need a new method developed thanks to a new approach to minimal 
surfaces, proposed by Ganchev \cite{G}. 

Here we explain with examples how this method works. Namely, in section 4 we prove that the
above mentioned families of polynomial surfaces of degree 6 from \cite{X-W-1} 
contain many congruent surfaces and that all the surfaces in these families 
coincide up to homothety and position in  space with the classical Enneper 
surface. We also show  how one can look for an explicit relation between  congruent 
minimal surfaces. Note that the construction of these examples in \cite{X-W-1} is
based on a theorem for the coefficients of a minimal polynomial surface of degree
six, but two of the equations in this theorem as formulated in \cite{X-W-1} are not exact. 
In section 3, we 
give the correct version of the theorem. Actually, the construction of the
above mentioned examples is not affected by the wrong equations.  
Unfortunately in \cite{X-W-1}  another general family is obtained by using 
one of the wrong equations and the surfaces in it are not minimal.
So in section 5 we construct a new large family of polynomial minimal surfaces
of degree 6 and we prove some relations between its surfaces.

\setcounter{equation}{0}
\section{Preliminaries}

Let $S$ be a regular surface defined by the parametric equation
$$
	\boldsymbol{ x=x}(u,v)=(x_1(u,v),x_2(u,v),x_3(u,v)) , \qquad 
	(u,v) \in U \subset \mathbb{R}^2 .
$$
Denote the derivatives of the vector function $\boldsymbol{x=x}(u,v)$ by
$	{\boldsymbol x}_u$, $	{\boldsymbol x}_v$, ${\boldsymbol x}_{uu}$, ${\boldsymbol x}_{uv}$, ${\boldsymbol x}_{vv}$.
The unit normal to the surface is defined by
$$
	{\boldsymbol U} = \frac{{\boldsymbol x}_u\times {\boldsymbol x}_v}{|{\boldsymbol x}_u\times {\boldsymbol x}_v|} .
$$
Then the coefficients of the first and the second fundamental form are given resp. by 
$$
	E={\boldsymbol x}_u^2, \qquad F={\boldsymbol x}_u{\boldsymbol x}_v , \qquad G={\boldsymbol x}_v^2 ,
$$
and
$$
	L={\boldsymbol U\,x}_{uu} , \qquad 	M={\boldsymbol U\,x}_{uv} , \qquad 	N={\boldsymbol U\,x}_{vv} .
$$

Recall that the parameters of the surface are called isothermal, if $E=G$, $F=0$.

The Gauss curvature $K$ and the mean curvature $H$ of a surface $S$ are 
defined by
$$
  K=\frac{LN-M^2}{EG-F^2} , \qquad	
  H=\frac{E N-2F M+G L}{2(EG-F^2)} .
$$
The surface $S$ is called {\it minimal} if its mean 
curvature vanishes identically. In this case the Gauss curvature is 
negative. The {\it normal curvature} $\nu$ of a minimal surface $S$ is defined 
to be the function $\nu=\sqrt{-K}$ , see \cite{G}. This is exactly the positive 
principal curvature of  $S$.

Let us now consider a minimal surface defined in isothermal parameters.
Then it can be considered as a real part of a complex curve. More precisely,
let $f(z)$ and $g(z)$ be two holomorphic functions. Define the Weierstrass 
complex curve $\Psi(z)$  by
\begin{equation} \label{eq:2.1}
   {\boldsymbol	\Psi}(z)=\int_{z_0}^z\left( \frac12 f(z)(1-g^2(z)),\frac{i}2 f(z)(1+g^2(z)), f(z)g(z) \right)dz \ . 
\end{equation}
Then $\boldsymbol\Psi(z)$ is a minimal curve, i.e. $\boldsymbol\Psi'^{\, 2}(z)=0$, and its real
and imaginary parts ${\boldsymbol x}(u,v)$ and ${\boldsymbol y}(u,v)$ define two  
minimal surfaces in isothermal parametrizations. We say that these two minimal 
surfaces are {\it conjugate}. Moreover, every minimal surface can be 
obtained at least locally as the real part of  a Weierstrass minimal curve.

\vspace{0.2cm}
{\bf Example.} Taking $f(z)=e^z$, $g(z)=e^{-z}$ we obtain a Weierstrass minimal curve
whose real part is the catenoid
$$
	{\boldsymbol x}(u,v)=( \cosh u\cos v, -\cosh u \sin v, u)
$$
and the imaginary part is the helicoid
$$
	{\boldsymbol y}(u,v)=(\sinh u\sin v , \sinh u\cos v , v) \ .
$$

\vspace{0.2cm}
Given two conjugate minimal surfaces ${\boldsymbol x}(u,v)$ and ${\boldsymbol y}(u,v)$ we can define the associated family
$$
	S_t \ \ :\ \qquad \boldsymbol{ assoc}_t(u,v)={\boldsymbol x}(u,v)\cos t+{\boldsymbol y}(u,v)\sin t \ .
$$
Then for any real number $t$ the surface $S_t$ is also minimal and
has the same first fundamental form as $S$.

It is easy to verify that the coefficients  of the first fundamental 
form of the associated family defined in the above way via the Weierstrass formula with functions 
$f(z)$, $g(z)$ are given by
\begin{equation} \label{eq:2.2}
	E = \frac14 | f | ^2 (1 + | g | ^2)^2 \qquad F=0 \qquad G = \frac14 | f | ^2 (1 + | g | ^2)^2 
\end{equation}
and the normal curvature is
\begin{equation} \label{eq:2.3}
	\nu=\frac{4|g'|}{|f|(1+|g|^2)^2} ,
\end{equation}
see \cite[Theorem~22.33]{G-A-S}.

A new approach to minimal surfaces was proposed by Ganchev in \cite{G}. 
He introduced special principal parameters and called them 
{\it canonical principal parameters}. The surface being  parametrized with 
them, the coefficients of its fundamental forms are given by
\begin{gather*}
	E=\frac1\nu \ , \qquad F=0 \ , \qquad G=\frac1\nu \\
	L=1 \ , \qquad M=0 \ , \qquad N=-1 \ .
\end{gather*}
Moreover, the surface in canonical principal parametrization is the real part of
a Weierstrass minimal curve generated by some 
functions $f(z)$, $g(z)$ with $f(z)=-1/g'(z)$, i.e. it is the real part of the 
special Weierstrass curve
$$
	\boldsymbol\Phi(z)=-\int_{z_0}^z  
	\left( \frac{1-g^2(z)}{2g'(z)}, \frac{i(1+g^2(z))}{2g'(z)}, \frac{g(z)}{g'(z)} \right) dz
$$

The following theorem is of great importance:

\begin{thm}\label{T:2.1} {\cite{G}}  
If a surface is parametrized with canonical principal parameters, 
then its normal curvature satisfies the differential equation
$$
	\Delta \ln \nu+2\nu=0 .
$$
Conversely, for any solution $\nu(u,v)$ of this equation (with $\nu_u\nu_v\ne0$), there exists
a {\bf unique} (up to position in space) minimal surface with normal 
curvature $\nu(u,v)$, $(u,v)$ being canonical principal parameters.
Moreover, the canonical principal parameters $(u,v)$ are determined uniquely up to 
the following transformations 
\begin{equation} \label{eq:2.4}
	\begin{array}{l}
		u=\varepsilon\bar u+a , \\
		v=\varepsilon\bar v+b ,
	\end{array} \qquad \varepsilon=\pm 1 \ , \ a=const., \ b=const. 
\end{equation}
\end{thm}

We will also use the following results:

\vspace{0.2cm}
\begin{thm}\label{T:2.2} \cite{OK}
Let the minimal surface $S$ be defined by the
real part of the Weierstrass minimal curve (\ref{eq:2.1}). Any solution 
of the differential equation 
\begin{equation} \label{eq:2.5}
	(z'(w))^2=-\frac1{f(z(w))g'(z(w))} \ .  
\end{equation}
defines a transformation of the isothermal parameters of $S$ to canonical principal 
parameters. Moreover the function $\tilde g(w)$ that defines $S$ via the Ganchev 
formula is given by $\tilde g(w)=g(z(w))$.
\end{thm}

\begin{thm}\label{T:2.5} \cite{OK} Let the holomorphic function $ g(z)$ generate a minimal 
surface in canonical principal parameters, i.e. via the Ganchev formula. 
Then, for an arbitrary complex number $\alpha$, and for an arbitrary real number 
$\varphi$, any of the functions
\begin{equation} \label{eq:2.6}
  e^{i\varphi}\frac{\alpha+g(z)}{ 1-\bar \alpha g(z)}   
\end{equation} 
generates the same surface  in canonical principal parameters (up to position 
in  space). Conversely,  any function that generates this surface (up to position 
in  space) in canonical principal parameters has the above form.
\end{thm}

More precisely, from the proof of the last theorem we see
that, up to a translation, the surface generated by the function (\ref{eq:2.6}) is
obtained from the surface generated by $g(z)$ after a rotation with the
matrix $AB$, where
 
$$
	A=\left(\begin{array}{ccc}
	     \cos\varphi  &  -\sin\varphi  &  0  \\
	     \sin\varphi  &   \cos\varphi  &  0  \\
	          0       &       0        &  1
	  \end{array}\right)  \ , \qquad  \qquad
	  	B=\left(\begin{array}{ccc}
	     \frac{1-a^2+b^2}{1+a^2+b^2}  &  -\frac{2ab}{1+a^2+b^2}  &  -\frac{2a}{1+a^2+b^2}  \vspace{.1cm} \\ 
	    -\frac{2ab}{1+a^2+b^2}  &  \frac{1+a^2-b^2}{1+a^2+b^2}   & -\frac{2b}{1+a^2+b^2}   \vspace{.1cm} \\
	     \frac{2a}{1+a^2+b^2}   &  \frac{2b}{1+a^2+b^2}          &  \frac{1-a^2-b^2}{1+a^2+b^2}
	  \end{array}\right)
$$
and $\alpha = a+ib$. In particular, if $\alpha=0$ this is a rotation by the angle $\varphi$
about the third axis.

\begin{thm}\label{T:2.6} \cite{OK} For arbitrary nonzero complex numbers $a$, $b$, the minimal surfaces
generated via the Weierstrass formula by the pairs of functions 
$$
	f(z)=a\,z^k   \ , \qquad g(z)=b\,z^n \ , \quad a, b \in \mathbb C \backslash\{0\} , \quad k,n\in\mathbb N \,  .  
$$  
coincide up to position in space and homothety.
\end{thm}

\setcounter{equation}{0}
\section{Parametric minimal surfaces of degree six.}
 
It is shown in \cite{X-W-1} that if a parametric minimal surface of degree 6
is presented in isothermal parameters, then it has the form
\begin{equation} \label{eq:3.1}
	\begin{array}{rl}{\boldsymbol r}(u,v)\ = & \ \ {\boldsymbol a}(u^6-15u^4v^2+15u^2v^4-v^6)+{\boldsymbol b}(3u^5v-10u^3v^3+3u v^5)\\
	                               & + {\boldsymbol c}(u^5-10u^3v^2+5uv^4)+{\boldsymbol d}(v^5-10u^2v^3+5u^4v)  \\
	                               & +{\boldsymbol e}(u^4-6u^2v^2+v^4)+{\boldsymbol f} uv(u^2-v^2)+{\boldsymbol g} u(u^2-3v^2)+{\boldsymbol h} v(v^2-3u^2) \\
	                               & +{\boldsymbol i} (u^2-v^2)+{\boldsymbol j} uv+{\boldsymbol k}u+{\boldsymbol l}v+{\boldsymbol m}
	             \end{array}               
\end{equation} 
where $\boldsymbol {a,\, b,\, c,\, d,\, e,\, f,\, g,\, h,\, i,\, j,\, k,\, l,\, m}$ are coefficient vectors. 
Then the authors of  \cite{X-W-1} formulate and give a sketch of the proof of a theorem about 
the relation between these  coefficients so that the surface be minimal. Actually, two of 
the equations in this theorem are not exact. Namely, the statement must look as follows

\begin{thm}\label{T:3.1.} The harmonic polynomial surface  (\ref{eq:3.1})    is a minimal  surface
if and only if its coefficient vectors satisfy the following system of equations
\begin{equation} \label{eq:3.2}
	\left\{\begin{array}{l}
			4\,\boldsymbol{ a}^2={\boldsymbol b}^2   \\
			{\boldsymbol a}.{\boldsymbol b}=0         \\
		  2\,{\boldsymbol a}.{\boldsymbol c}-{\boldsymbol b}.{\boldsymbol d}=0  \\
		  2\,{\boldsymbol a}.{\boldsymbol d}+{\boldsymbol b}.{\boldsymbol c}=0  \\ 
		  25\, {\boldsymbol c}^2-25\,{\boldsymbol d}^2+48\,{\boldsymbol a}.{\boldsymbol e}-6\,{\boldsymbol b}.{\boldsymbol f}=0 \\
		  25\, {\boldsymbol d}.{\boldsymbol c} +12\,{\boldsymbol b}.{\boldsymbol e}+6{\boldsymbol a}.{\boldsymbol f}=0          \\
		  16\,{\boldsymbol e}^2-{\boldsymbol f}^2+30\,{\boldsymbol c}.{\boldsymbol g}+30\,{\boldsymbol d}.{\boldsymbol h}+24\,{\boldsymbol a}.{\boldsymbol i}-6\, {\boldsymbol b}.{\boldsymbol j}=0  \\
		  4\,{\boldsymbol e}.{\boldsymbol f}-15\,{\boldsymbol c}.{\boldsymbol h}+15\,{\boldsymbol d}.{\boldsymbol g}+6\,{\boldsymbol b}.{\boldsymbol i}+6\,{\boldsymbol a}.{\boldsymbol j}=0   \\
		  9\,{\boldsymbol g}^2-9\,{\boldsymbol h}^2+16\,{\boldsymbol e}.{\boldsymbol i}-2\,{\boldsymbol f}.{\boldsymbol j}+10\,{\boldsymbol c}.{\boldsymbol k}-10\,{\boldsymbol l}.{\boldsymbol d}=0   \\
		  9\,{\boldsymbol g}.{\boldsymbol h}-2\,{\boldsymbol f}.{\boldsymbol i}-4\,{\boldsymbol e}.{\boldsymbol j}-5\,{\boldsymbol d}.{\boldsymbol k}-5\,{\boldsymbol c}.{\boldsymbol l}=0   \\
		  4\,{\boldsymbol i}^2-{\boldsymbol j}^2+6\,{\boldsymbol g}.{\boldsymbol k}+6\,{\boldsymbol h}.{\boldsymbol l}=0    \\
		  2\,{\boldsymbol i}.{\boldsymbol j}+3\,{\boldsymbol g}.{\boldsymbol l}-3\,{\boldsymbol h}.{\boldsymbol k}=0        \\
		  18\,{\boldsymbol a}.{\boldsymbol g}+9\,{\boldsymbol b}.{\boldsymbol h}+20\,{\boldsymbol e}.{\boldsymbol c}-5\,{\boldsymbol f}.{\boldsymbol d}=0   \\
		  18\,{\boldsymbol a}.{\boldsymbol h}-9\,{\boldsymbol b}.{\boldsymbol g}-20\,{\boldsymbol e}.{\boldsymbol d}-5\,{\boldsymbol f}.{\boldsymbol c}=0   \\
		  6\,{\boldsymbol a}.{\boldsymbol k}-3\,{\boldsymbol b}.{\boldsymbol l}+10\,{\boldsymbol c}.{\boldsymbol i}-5\,{\boldsymbol d}.{\boldsymbol j}+12\,{\boldsymbol e}.{\boldsymbol g}+3\,{\boldsymbol f}.{\boldsymbol h}=0   \\
		  6\,{\boldsymbol a}.{\boldsymbol l}+3\,{\boldsymbol b}.{\boldsymbol k}+5\,{\boldsymbol c}.{\boldsymbol j}+10\,{\boldsymbol d}.{\boldsymbol i}+3\,{\boldsymbol f}.{\boldsymbol g}-12\,{\boldsymbol e}.{\boldsymbol h}=0   \\
		  4\,{\boldsymbol e}.{\boldsymbol k}-{\boldsymbol f}.{\boldsymbol l}+3\,{\boldsymbol h}.{\boldsymbol j}+6\,{\boldsymbol g}.{\boldsymbol i}=0   \\
		  4\,{\boldsymbol e}.{\boldsymbol l}+{\boldsymbol f}.{\boldsymbol k}+3\,{\boldsymbol g}.{\boldsymbol j}-6\,{\boldsymbol h}.{\boldsymbol i}=0   \\
		  2\,{\boldsymbol l}.{\boldsymbol i}+{\boldsymbol k}.{\boldsymbol j}=0    \\
		  2\,{\boldsymbol k}.{\boldsymbol i}-{\boldsymbol l}.{\boldsymbol j}=0    \\
		  {\boldsymbol k}^2={\boldsymbol l}^2    \\
		  {\boldsymbol k}.{\boldsymbol l}=0
	\end{array}  \right.
\end{equation}  
\end{thm}

The proof is based on the vanishing of the coefficients of 
$u^\alpha v^\beta$ in the power expansion of the isothermal condition $E=G$, $F=0$. 
The differences with the theorem as formulated in \cite{X-W-1} are in the seventh and the
twelfth equations.

Using the first four equations of the system (\ref{eq:3.2}), we can suppose that 
the coefficient vectors are ${\boldsymbol a}=(a_1,a_2,0)$, ${\boldsymbol b}=(-2a_2,2a_1,0)$, 
${\boldsymbol c}=(c_1,c_2,c_3)$, ${\boldsymbol d}=(-c_2,c_1,d_3)$, ${\boldsymbol e}=(e_1,e_2,e_3)$,
${\boldsymbol f}=(f_1,f_2,f_3)$,  ${\boldsymbol g}=(g_1,g_2,g_3)$, ${\boldsymbol h}=(h_1,h_2,h_3)$, 
${\boldsymbol i}=(i_1,i_2,i_3)$, ${\boldsymbol j}=(j_1,j_2,j_3)$, ${\boldsymbol k}=(k_1,k_2,k_3)$, 
${\boldsymbol l}=(l_1,l_2,l_3)$, ${\boldsymbol m}=(0,0,0)$,  see \cite{X-W-1}.

\setcounter{equation}{0}
\section{First example.}

If we suppose $\boldsymbol{ c=d=g=h=k=l=o}$, $j_1=2i_2$, $j_2=-2i_1$, we
obtain a family of minimal surfaces, depending
on four shape parameters $a_1,a_2,i_1,i_2$. Two special cases
are studied in \cite{X-W-1}. Namely, if $a_2=i_2=0$ the surface is 
$$
	{\boldsymbol r_1}[a_1,i_1](u,v)=\big( X_1(u,v),Y_1(u,v),Z_1(u,v)\big)
$$
where
$$
  \begin{array}{l}
	     X_1(u,v)=a_1 (u^6 - 15 u^4 v^2 + 15 u^2 v^4-v^6) + i_1 (u^2 -  v^2)  \vspace{2mm} \\
		   Y_1(u,v)=2a_1(3u^5v-10u^3v^3+3u v^5)-2i_1 u v   \vspace{2mm} \\ 
		   Z_1(u,v)=\sqrt{\frac32}\,\sqrt{|a_1i_1|-a_1i_1\,}\ (u^4-6u^2v^2+v^4)-2\sqrt 6\,
		            \sqrt{|a_1i_1|+a_1i_1\,}\ u v(u^2-v^2)
  \end{array}
$$  
and if $a_1=i_1=0$, the surface is 
$$
	{\boldsymbol r_2}[a_2,i_2](u,v)=\big( X_2(u,v),Y_2(u,v),Z_2(u,v)\big)
$$
where
$$
  \begin{array}{l}
	    X_2(u,v)=-2a_2(3u^5v-10u^3v^3+3u v^5)+2i_2 u v   \vspace{2mm} \\ 
	    Y_2(u,v)=a_2 (u^6 - 15 u^4 v^2 + 15 u^2 v^4-v^6) + i_2  (u^2 -  v^2)  \vspace{2mm} \\ 
		  Z_2(u,v)= \sqrt{\frac32}\,\sqrt{|a_2i_2|-a_2i_2\, }\ (u^4-6u^2v^2+v^4)-2\sqrt 6\,
		            \sqrt{|a_2i_2|+a_2i_2\, }\ u v(u^2-v^2)
  \end{array} 
$$  

The following theorem is proved in \cite{X-W-1}:

\begin{thm}\label{T:4.1.} 
{\it The surfaces defined by the charts ${\boldsymbol r_1}[a_1,i_1](u,v)$ and 
${\boldsymbol r_2}[-a_1,i_1](u,v)$  are conjugate minimal surfaces.}
\end{thm}

The figures bellow represent the surfaces ${\boldsymbol r_1}[1,500](u,v)$ and ${\boldsymbol r_2}[-1,500](u,v)$ for $u,v\in [-4,4]$.

\begin{tabular}{cc}
		\qquad\qquad\epsfig{file=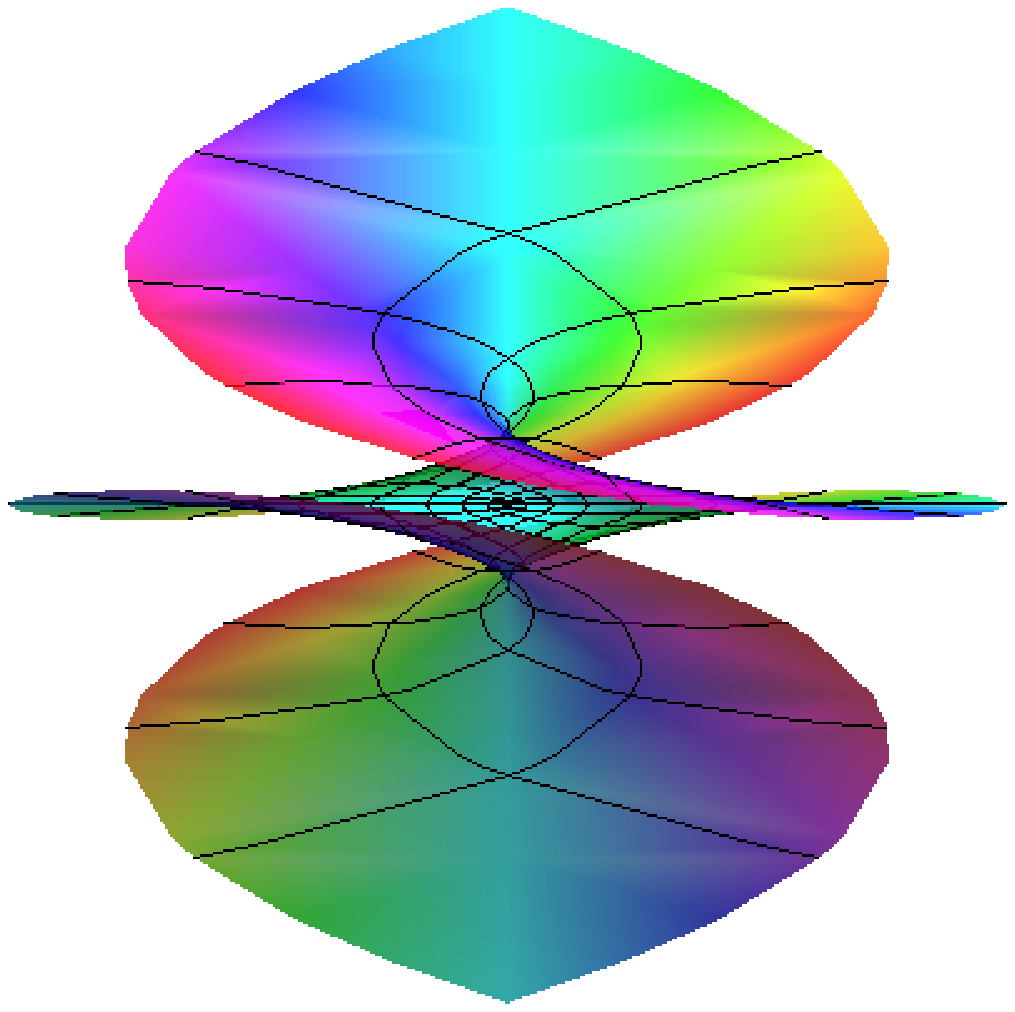,width=0.28\textwidth}   &
    \hspace{2cm}  \qquad \epsfig{file=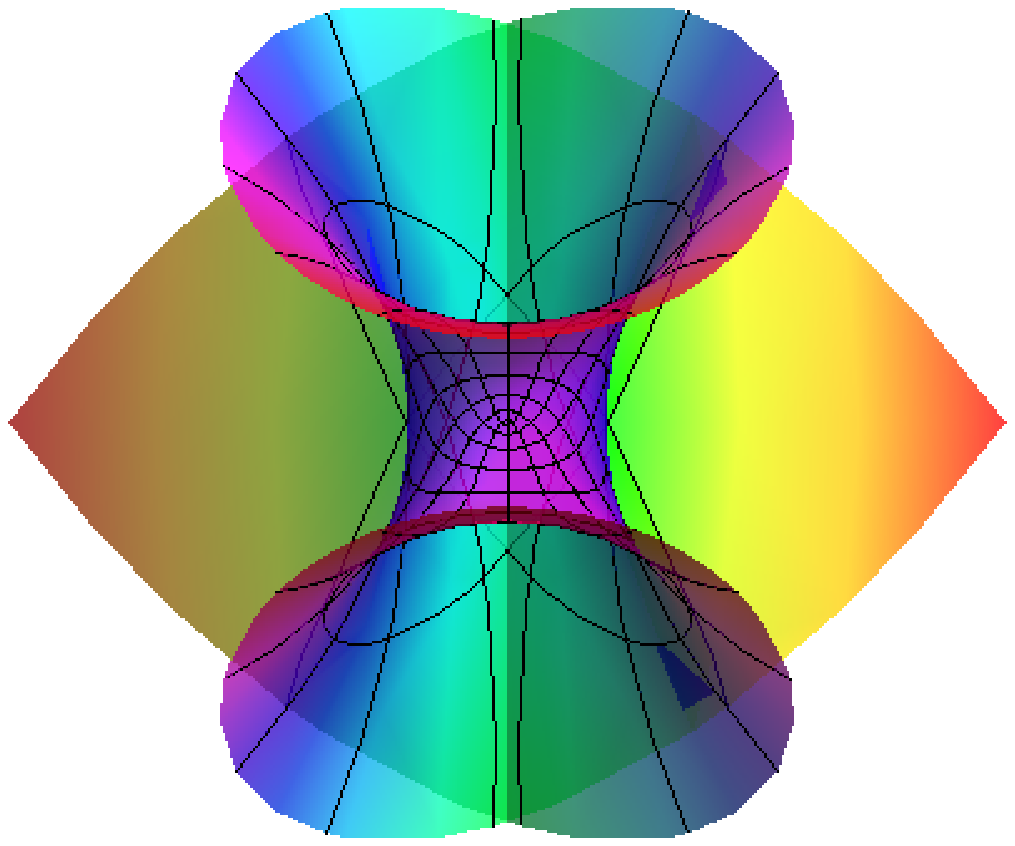,width=0.28\textwidth} \\ 
    \qquad\quad{\bf Fig. 4.1}: ${\boldsymbol r_1}[1,500](u,v)$ & \hspace{2cm}   {\bf Fig. 4.2}: ${\boldsymbol r_1}[1,-500](u,v)$ \\
\end{tabular}

\vspace{0.2cm}

It turns out that the conjugate surfaces in Theorem \ref{T:4.1.} actually coincide. We will see 
that moreover they coincide with any surface of their associated family. Namely, we have

\begin{thm}\label{T:4.2.} For any real number \ $t$ \ the surface 
\begin{equation} \label{eq:4.6}
	S_t \ \ :\ \qquad \boldsymbol{ assoc}_t[a_1,i_1](u,v)={\boldsymbol r_1}[a_1,i_1]\cos t+{\boldsymbol r_2}[-a_1,i_1](u,v)\sin t
\end{equation}  
coincides up to position in  space with that defined by \ ${\boldsymbol r_1}[a_1,i_1](u,v)=\boldsymbol{ assoc}_0[a_1,i_1](u,v)$.
\end{thm}

{\it Proof.}
Suppose for example that $a_1i_1>0$. We will change the parameters of $S_t$ 
to canonical principal ones. First we need the functions that generate 
$\boldsymbol{ assoc}_t[a_1,i_1](u,v)$ via the Weierstrass formula. They can be found 
from the well known relations
\begin{equation} \label{eq:4.8}
	f_t=\varphi_1-i\varphi_2\ , \qquad  g_t=\frac{\varphi_3}{\varphi_1-i\varphi_2} \ ,
\end{equation}
where 
$$
	(\varphi_1,\varphi_2,\varphi_3)= (\boldsymbol{ assoc}_t)_u-i\,(\boldsymbol{ assoc}_t)_v \ ,
$$
see e.g. \cite{G-A-S}, Lemma 22.31. In our case these equalities give
\begin{equation} \label{eq:4.1}
	f_t(z)=4i_1e^{-it}z  \ , \qquad g_t(z)=\frac{\sqrt 3\, a_1 i}{\sqrt{a_1i_1}}z^2 \ .
\end{equation} 
With these functions a solution of (\ref{eq:2.5}) is 
$$
	z(w)=	\frac{\sqrt{(1+i)e^{it/2}w}}{\sqrt{2}\, (3a_1i_1)^{1/8}} \ .
$$
So we change $z$ in $g_t(z)$ by $z(w)$  and we obtain that the function
$$
	\tilde g_t(w)=\frac{(-1+i)3^{1/4}a_1e^{i t/2}w}{2(a_1i_1)^{3/4}}
$$
generates $S_t$ in canonical principal parameters. Then, using (\ref{eq:2.3}) 
we obtain the normal curvature in canonical principal parameters to be
(recall that $t$ is real) 
$$
	\tilde \nu=\frac{4|\tilde g_t'|^2}{(1+|\tilde g_t|^2)^2}
	=8\sqrt{3}\frac{\left|\frac{a_1}{(a_1i_1)^{3/4}}\right|^2}{\left(2+\sqrt{3}\left|\frac{a_1}{(a_1i_1)^{3/4}}w\right|^2\right)^2} 
$$
Since this number does not depend on $t$, Theorem \ref{T:2.1} implies that any 
surface $S_t$ is congruent with every other surface in the associated family and 
in particular to $S_0$.
\hfill{\qed}

\vspace{0.1cm}
In the same way we can prove

\begin{thm}\label{T:4.3.} 
{\it The minimal surfaces defined by the charts ${\boldsymbol r_1}[a_1,i_1](u,v)$,
${\boldsymbol r_1}[a_1,-i_1](u,v)$, ${\boldsymbol r_1}[-a_1,i_1](u,v)$, ${\boldsymbol r_1}[-a_1,-i_1](u,v)$ are congruent.}
\end{thm}

\begin{cor}\label{T:4.4.} All the minimal surfaces from the associated families
$ \boldsymbol{ assoc}_t[a_1,i_1](u,v) $, $ \boldsymbol{ assoc}_t[a_1,-i_1](u,v) $, $ \boldsymbol{ assoc}_t[-a_1,i_1](u,v) $,
$ \boldsymbol{ assoc}_t[-a_1,-i_1](u,v) $ are congruent.
\end{cor}

Of course, the proof of Theorem \ref{T:4.2.} given above hides the relation between the charts 
${\boldsymbol r_1}[a_1,i_1](u,v)$ and $\boldsymbol{ assoc}_t[a_1,i_1](u,v)$. To show this relation we 
first remark that we can find an alternative end of the proof of Theorem \ref{T:4.2.}  if we notice 
that the generating functions $\tilde g_0(z)$ and $\tilde g_t(z)$ satisfy 
\begin{equation} \label{eq:4.2}
	\tilde g_t(z)=e^{\frac{it}2}g_0(z) \ .
\end{equation} 
Indeed,   the last equality and Theorem \ref{T:2.5} imply the conclusion of Theorem \ref{T:4.2.}.
Moreover in view of (\ref{eq:4.2}) and Theorem \ref{T:2.5} we can say that up to a translation  
$S_t$ is obtained from $S_0$ with a rotation of angle $t/2$ around the third axis. To realize 
this rotation we first obtain the Weierstrass minimal curves of   ${\boldsymbol r_1}[a_1,i_1](u,v)$
and $\boldsymbol{ assoc}_t[a_1,i_1](u,v)$. By  (\ref{eq:2.1}), (\ref{eq:4.1}), and the definitions of 
 ${\boldsymbol r_1}[a_1,i_1](u,v)$ and ${\boldsymbol r_2}[-a_1,i_1](u,v)$,  these curves are
\begin{equation} \label{eq:4.3}
	\left(z^2 (i_1 + a_1 z^4), i z^2 (i_1 - a_1 z^4), i  \sqrt{3a_1 i_1}\, z^4\right)
\end{equation} 
\begin{equation} \label{eq:4.4}
	\left(e^{-i t} z^2 (i_1 + a_1 z^4), i e^{-i t} z^2 (i_1 - a_1 z^4), 
   i  e^{-i t}\, \sqrt{3a_1 i_1} z^4\right)
\end{equation} 
respectively. So we see that the transformation from the chart (\ref{eq:4.3}) to 
the chart (\ref{eq:4.4}) can be made without any translation, but we have to make  
the change of the parameters $z\ \rightarrow \ e^{-i t/4} z $ in (\ref{eq:4.3}).
So we obtain the following parametrization of the minimal curve defined by (\ref{eq:4.3}):
$$
	\left(e^{-i t/2} z^2 (i_1 + a_1 e^{-i t} z^4), 
 i e^{-i t/2} z^2 (i_1 - a_1 e^{-i t} z^4), 
 i e^{-i t} \sqrt{3a_1 i_1} z^4\right)
$$
Now the rotation of the last chart with angle $t/2$ around the third axis gives (\ref{eq:4.4}). 
The real form of the above change of parameters is
\begin{equation} \label{eq:4.5}
	u\ \rightarrow \ u \cos\frac t4 + v \sin\frac t4    \qquad v\ \rightarrow \ -u \sin\frac t4 + v\cos\frac t4 \ .
\end{equation} 
As a result we obtain
\begin{thm}\label{T:4.5.} 
{\it For any real number $t$ the chart $\boldsymbol{ assoc}_t[a_1,i_1](u,v)$
is obtained from ${\boldsymbol r_1}[a_1,i_1](u,v)$ with a change of the parameters (\ref{eq:4.5})
and then a rotation with angle $t/2$ aroud the third axis.}
\end{thm} 

Remark also that from (\ref{eq:4.1}) and Theorem \ref{T:2.6}
we obtain

\begin{thm}\label{T:4.6.} 
{\it For any real number \ $t$ \ and for arbitrary values of the shape parameters 
$a_1$, $i_1$,  the surface (\ref{eq:4.6}) coincides up to position in  space and 
homothety with the surface generated via the Weierstrass formula by the functions
\begin{equation} \label{eq:4.7}
	f(z)=z  \ , \qquad g(z)=z^2 \ .   
\end{equation}}
\end{thm} 

In particular, Theorem \ref{T:4.6.} says that all the surfaces  in the families 
${\boldsymbol r_1}[a_1,i_1](u,v)$ and $ {\boldsymbol r_2}[a_2,i_2](u,v)$ coincide up to position 
in  space and homothety.

At the end of this section we note that, for any values of the shape parameters, the
function $\tilde g(z)$ that generates the considered surfaces has the form $\tilde g(z)=Cz$ with
a complex constant $C$, just like the surfaces homothetic to the Enneper surface. This looks strange for a 
polynomial surface of degree 6. Now we are obliged to take a more careful look on the equations
of  ${\boldsymbol r}_1[a_1,i_1](u,v)$ and ${\boldsymbol r}_2[a_2,i_2](u,v)$. We remark  that they are all homothetic
to the Enneper surface. If for example $a_1i_1>0$, to realize the correspondence it suffices to make  
the change of the parameters
$$
	u^2-v^2 \rightarrow \sqrt{\frac{i_1}{a_1}}\,U  \qquad\qquad 2uv   \rightarrow \sqrt{\frac{i_1}{a_1}}\,V \ .
$$

\setcounter{equation}{0}
\section{Second example.}

In this section we suppose that $\boldsymbol {g=h=i=j=k=l=o}$. Then, with the assumptions 
from the end of section 3, the system (\ref{eq:3.2}) is equivalent to
\begin{equation} \label{eq:5.1}
	\left\{\begin{array}{l} 
		  25\, {\boldsymbol c}^2-25\,{\boldsymbol d}^2+48\,{\boldsymbol a}.{\boldsymbol e}-6\,{\boldsymbol b}.{\boldsymbol f}=0 \\
		  25\, {\boldsymbol c}.{\boldsymbol d} +12\,{\boldsymbol b}.{\boldsymbol e}+6{\boldsymbol a}.{\boldsymbol f}=0          \\
		  16\,{\boldsymbol e}^2-{\boldsymbol f}^2=0  \\
		  \,{\boldsymbol e}.{\boldsymbol f}=0   \\
		  4\,{\boldsymbol c}.{\boldsymbol e}-\,{\boldsymbol d}.{\boldsymbol f}=0   \\
		  4\,{\boldsymbol d}.{\boldsymbol e}+\,{\boldsymbol c}.{\boldsymbol f}=0   \\
	\end{array}  \right.
\end{equation}  
Suppose in addition that the vectors $\boldsymbol e$ and $\boldsymbol f$ lie in $span\{\boldsymbol {a, b}\}$.
If the surface is not  part of a plane, then, using $(\ref{eq:5.1})_3-(\ref{eq:5.1})_6$,
we find ${\boldsymbol e}=(e_1,e_2,0)$, ${\boldsymbol f}=(4e_2,-4e_1,0)$. So (\ref{eq:5.1}) takes the form
\begin{equation} \label{eq:5.2}
	\left\{\begin{array}{l} 
		  25(c_3^2-d_3^2)+96(a_1e_1+a_2e_2)=0 \\
		  25c_3d_3+48(a_1e_2-a_2e_1)=0          \\
		  c_1e_1+c_2e_2=0   \\
		  c_1e_2-c_2e_1=0   \\
	\end{array}  \right.
\end{equation}  
If $e_1=e_2=0$, then the surface is part of a plane. So we assume
that $(e_1,e_2)\ne (0,0)$. Then, from  (\ref{eq:5.2})  we obtain
$c_1=c_2=0$,
$$
	e_1=\frac{25(-a_1c_3^2+2a_2c_3d_3+a_1d_3^2)}{96(a_1^2+a_2^2)}  \qquad
	e_2=\frac{25(-a_2c_3^2-2a_1c_3d_3+a_2d_3^2)}{96(a_1^2+a_2^2)}
$$
So, as in section 4, we obtain a family of minimal surfaces depending on the four
shape parameters $a_1,a_2,c_3,d_3$:
$$
	{\boldsymbol s}[a_1,a_2,c_3,d_3](u,v)=\big( X(u,v),Y(u,v),Z(u,v)\big)
$$
where
$$
  \begin{array}{rcl}
	     X(u,v)& = &a_1 (u^6 - 15 u^4 v^2 + 15 u^2 v^4-v^6)-2a_2\,uv(3u^4-10u^2v^2+3v^4)   \vspace{2mm} \\   
	    \ds    & &-\frac{25(a_1c_3^2-2a_2c_3d_3-a_1d_3^3)}{96(a_1^2+a_2^2)}\,(u^4-6u^2v^2+v^4) 
	    -\frac{25(a_2c_3^2+2a_1c^3d_3-a_2d_3^2)}{24(a_1^2+a_2^2)} \, uv(u^2-v^2)    
	    \vspace{2mm} \\ 	     	     	     
		   Y(u,v) & = &2a_1\,uv(3u^4-10u^2v^2+3 v^4)+a_2\,(u^6-15u^4v^2+15u^2v^4-v^6)  \\
		 \displaystyle    & & +\frac{25(a_1c_3^2-2a_2c_3d_3-a_1d_3^3)}{24(a_1^2+a_2^2)}\, uv(u^2-v^2) 
		  -\frac{25(a_1c_3^2-2a_2c_3d_3-a_1d_3^3)}{96(a_1^2+a_2^2)}\,(u^4-6u^2v^2+v^4)  
		  \vspace{2mm} \\ 		   
		   Z(u,v) & = & c_3\,u(u^4-10u^2v^2+5v^4)+d_3\,v(5u^4-10u^2v^2+v^4) \ ,
  \end{array}
$$  

Now we shall consider the surfaces defined by
$$
	{\boldsymbol s}_1[a_1,c_3]={\boldsymbol s}[a_1,0,c_3,0](u,v)
$$
and
$$
	{\boldsymbol s_2}[a_2,d_3]={\boldsymbol s}[0,a_2,0,d_3](u,v)
$$
Using (\ref{eq:4.8}) we obtain the functions that generate these surfaces via the Weierstrass
formula. They are respectively
\begin{equation} \label{eq:5.3}
	f_1(z)=-\frac{25 c_3^2}{12a_1}\, z^3    \qquad
	g_1(z)=-\frac{12a_1}{5c_3}\, z  \ ,
\end{equation}
\begin{equation} \label{eq:5.4}
	f_2(z)=-\frac{25 d_3^2}{12a_2}\, i z^3    \qquad
	g_2(z)=\frac{12a_2}{5d_3}\, z  \ .
\end{equation}
From  (\ref{eq:2.1}) and (\ref{eq:5.3}) we find the Weierstrass minimal curve,
whose real part is ${\boldsymbol s}_1[a_1,c_3]$:
$$
	\Psi_1(z)=\left(  \big(a_1 z^2-\frac{25 c_3^2 }{96 a_1}\big)z^4, 
	-   i\big( a_1 z^2+\frac{25 c_3^2 }{96 a_1}\big)z^4 , c_3 z^5\right)
$$
Then it is easy to check that the imaginary part of $\Psi_1(z)$ is
${\boldsymbol s}_2[-a_1,c_3]$. Hence we have

\begin{thm}\label{T:5.1.} 
{\it The minimal surfaces defined by the charts ${\boldsymbol s_1}[a_1,c_3](u,v)$ and
${\boldsymbol s_2}[-a_1,c_3](u,v)$ are conjugate.}
\end{thm}

Concerning the pairs (\ref{eq:5.3}) and (\ref{eq:5.4}), we take the following solutions 
of the system  (\ref{eq:2.5}) for the transition to canonical principal parameters:
$$
	z_1(w)=\left(\frac{i\sqrt{5}\,w}{2\sqrt{c_3}}\right)^{2/5}  \qquad
	z_2(w)=i\left(\frac{i\sqrt{5}\,w}{2\sqrt{d_3}}\right)^{2/5} \,
$$
respectively. Hence the functions
$$
	\tilde g_1(z)=-\frac{6.\,2^{3/5}a_1}{5^{4/5}c_3}\left(\frac{iw}{\sqrt{c_3}}\right)^{2/5}
	\qquad\qquad
	\tilde g_2(z)=\frac{6i.\,2^{3/5}a_2}{5^{4/5}d_3}\left(\frac{iw}{\sqrt{d_3}}\right)^{2/5}
$$
generate the surfaces in canonial principal parameters. Hence
it is clear that the normal curvatures of ${\boldsymbol s_1}[a_1,c_3](u,v)$ and
${\boldsymbol s_2}[-a_1,c_3](u,v)$  coincide, and so we have

\begin{thm}\label{T:5.2.} 
{\it The conjugate minimal surfaces defined by the charts ${\boldsymbol s_1}[a_2,c_3](u,v)$
and ${\boldsymbol s_2}[-a_1,c_3](u,v)$ are congruent.}
\end{thm}

\begin{cor}\label{T:5.3.} Any two minimal surfaces from the associated family
$$
	{\boldsymbol s_1}[a_1,c_3](u,v)\cos t+ {\boldsymbol s_2}[-a_1,c_3](u,v)\sin t  \qquad t\in [0,\pi/2]
$$
are congruent.
\end{cor}

To find the exact relation between the conjugate charts ${\boldsymbol s_1}[a_1,c_3](u,v)$ and
${\boldsymbol s_2}[-a_1,c_3](u,v)$, we can use the same arguments as in 
section 4. Then we can see that the chart ${\boldsymbol s_2}[-a_1,c_3](u,v)$ 
is obtained from ${\boldsymbol s_1}[a_1,c_3](u,v)$ with the change of the parameters
$u \rightarrow v$, $v \rightarrow -u$ and then a rotation with angle $\pi/2$
about the third axis.

Of course, this is not the only way to obtain ${\boldsymbol s_2}[-a_1,c_3](u,v)$ 
from ${\boldsymbol s_1}[a_1,c_3](u,v)$. For example, another possibility is
to make the following change of the parameters in ${\boldsymbol s_1}[a_1,c_3](u,v)$:
$$
	u \rightarrow  \sqrt{\frac{5+\sqrt 5}8}\,u +\frac{\sqrt 5-1}4\,v  
	\qquad\qquad
	v \rightarrow  -\frac{\sqrt 5-1}4\,u+\sqrt{\frac{5+\sqrt 5}8}\,v 
$$
and to rotate about the third axis with an angle $\varphi$, such that
$$
	\cos\varphi=\sqrt{\frac{5+\sqrt 5}8}
	\qquad\qquad
	\sin\varphi= \frac{-1+\sqrt 5}4
$$
An analogous remark holds for the surfaces in section 4.

Now consider the general chart ${\boldsymbol s}[a_1,a_2,c_3,d_3](u,v)$.
Using (\ref{eq:4.8}) we see that it is generated via the Weierstrass
formula with the functions 
\begin{equation} \label{eq:5.5}
	f(z)=-\frac{25(c_3-id_3)^2}{12(a_1+ia_2)}z^3   \qquad\quad
	g(z)=-\frac{12(a_1+ia_2)}{5(c_3-id_3)}z \ .
\end{equation}
A corresponding solution of the differential equation (\ref{eq:2.5})  is
$$
	z(w)=\frac{5^{1/5}(-w^2)^{1/5}}{2^{2/5}(c_3-id_3)^{1/5}}
$$
which gives the following generating function in canonical principal parameters:
$$
	\tilde g(w)=-\frac{6.2^{3/5}(a_1+ia_2)}{5^{4/5}(c_3-id_3)^{6/5}}(-w^2)^{1/5} 
$$
and the following normal curvature in these parameters:
$$
	\frac{1152.50^{1/5}\frac{a_1^2+a_2^2}{(c_3^2+d_3^2)^{6/5}|w|^{6/5}}}
	     {\left(25+72.50^{1/5}\frac{(a_1^2+a_2^2)|w|^{4/5}}{(c_3^2+d_3^2)^{6/5}}\right)^2}
$$
Hence it follows that every two surfaces for which
$$
	\frac{(a_1^2+a_2^2)^5}{(c_3^2+d_3^2)^6}
$$  
is one and the same number are congruent. In particular, this implies Theorem \ref{T:5.2.}.

On the other hand, using (\ref{eq:5.5}) and Theorem \ref{T:2.6}, we obtain

\begin{thm}\label{T:5.4.} 
{\it For any shape parameters, the surface defined by ${\boldsymbol s}[a_1,a_2,c_3,d_3](u,v)$ coincides with
the surface generated via the Weierstrass formula by the functions
$$	f(z)=z^3  \ , \qquad g(z)=z \ $$  
up to position in  space and homothety.}
\end{thm} 

Figure 5.1 represents such a surface from different points of view.

\begin{tabular}{c}
		\qquad\quad\epsfig{file=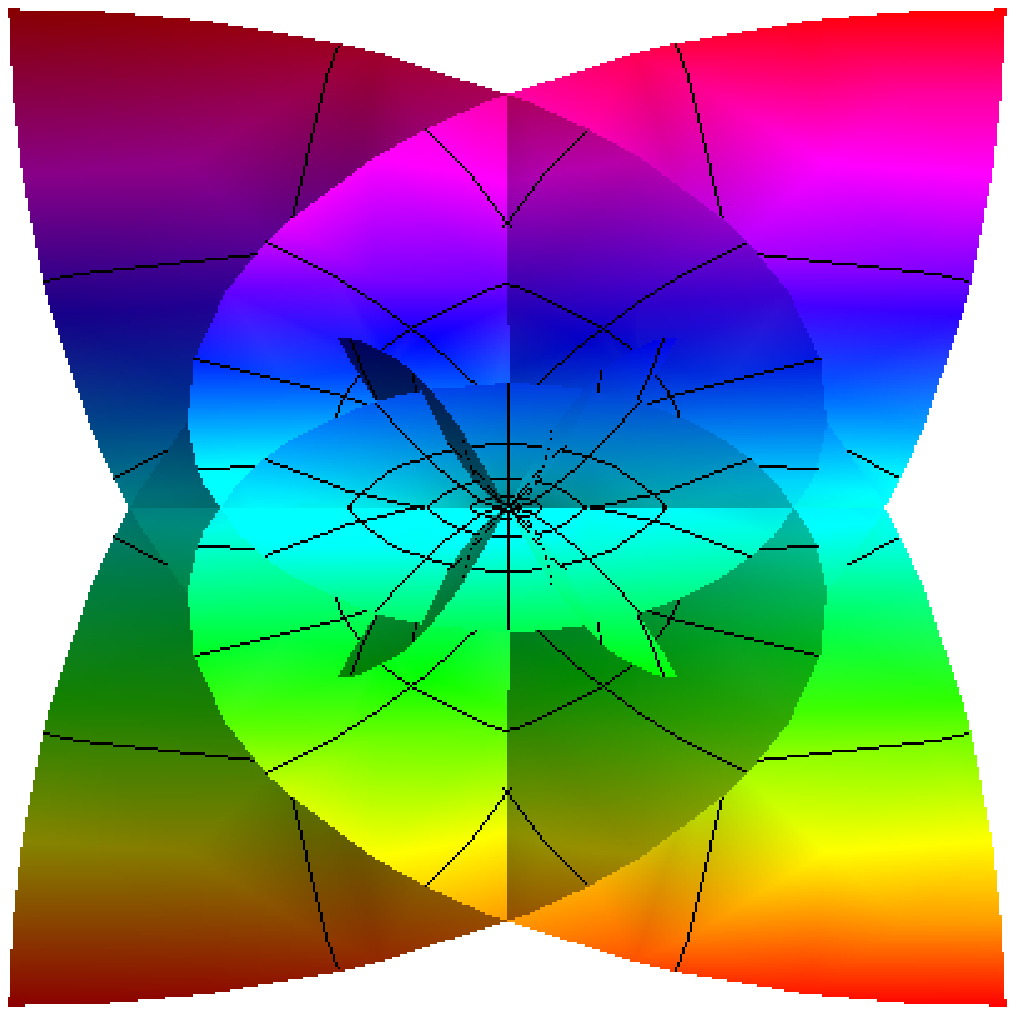,width=0.28\textwidth}   
    \hspace{2cm}  \qquad \epsfig{file=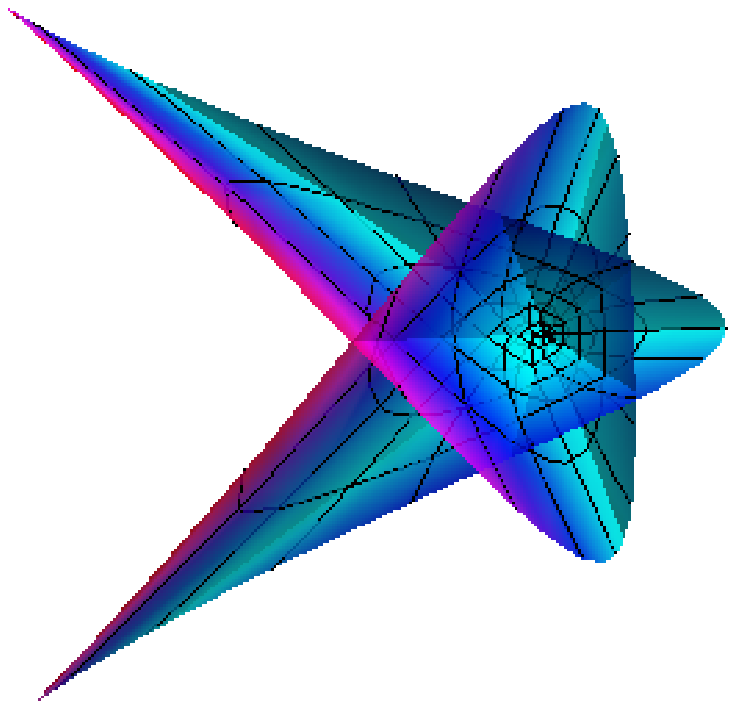,width=0.28\textwidth} \\ 
      {\bf Fig. 5.1}: ${\boldsymbol s_1}[3,0,1,0](u,v)$ for $u,v\in [-0.1,0.1]$. \\
\end{tabular} 

\vspace{0.2cm}
At the end of this article we note that the parametric minimal surfaces
of degree six in section 4 are completely different from those in section 5. 
This follows from their normal curvatures and Theorem \ref{T:2.1}, or 
from the functions that generate these surfaces in canonical principal 
parameters and Theorem \ref{T:2.5}. Indeed, any of the surfaces in section 4 
is generated in canonical principal parameters by
a function of the form $\tilde g(z)=C z$, whereas any surface defined by a chart
${\boldsymbol s}[a_1,a_2,c_3,d_3](u,v)$ is generated by a function of the
form $\tilde g(z)=C (-z^2)^{1/5}$.

\end{document}